\documentclass{article}

\makeatletter
\newcommand{\globalcolor}[1]{%
  \color{#1}\global\let\default@color\current@color
}
\makeatother

\usepackage{rotating}
\usepackage{arxiv}

\usepackage[utf8]{inputenc} 
\usepackage[T1]{fontenc}    
\usepackage{hyperref}       
\usepackage{url}            
\usepackage{booktabs}       
\usepackage{amsfonts}       
\usepackage{nicefrac}       
\usepackage{lipsum}

\usepackage{float}
\usepackage{morefloats}
\usepackage{color}
\usepackage{multirow}
\usepackage{latexsym}
\usepackage{amsmath,amssymb,amsthm,graphicx}
\usepackage{dsfont}
\usepackage{enumitem}
\usepackage{hyperref}

\newtheoremstyle{thm}
{9pt}
{9pt}
{\itshape}
{}
{\bfseries}
{.}
{ }
{}
\theoremstyle{thm}

\newtheorem{theorem}{Theorem}[section]

\newtheorem{corollary}[theorem]{Corollary}

\newcommand{\edist}{\stackrel{{\cal{D}}}{=}}
\newcommand{\e}{\text{e}}
\newcommand{\be}{\begin{eqnarray}}
\newcommand{\ee}{\end{eqnarray}}
\newcommand{\bewend}{\hspace*{2mm}\rule{3mm}{3mm}}
\newcommand{\verk}{\stackrel{{\cal D}}{\longrightarrow}}
\newcommand{\BE}{\mathbb{E}}

\newtheoremstyle{def}
{9pt}
{9pt}
{}
{}
{\bfseries}
{.}
{ }
{}
\theoremstyle{def}

\newtheorem{remark}[theorem]{Remark}

\newcommand{\E}{\mathbb{E}} 
\newcommand{\PP}{\mathbb{P}} 
\newcommand{\ii}{\textrm{i}}
\newcommand{\RR}{\mathbb{R}}
\newcommand{\re}{\textrm{Re}}
\newcommand{\ri}{\textrm{RI}}
\newcommand{\im}{\textrm{Im}}

\renewcommand{\footnoterule}{%
	\kern -3.5pt
	\hrule width \textwidth height 1pt
	\kern 3.5pt
}

\makeatletter
\def\blfootnote{\xdef\@thefnmark{}\@footnotetext}
\makeatother

\title{Tests for circular symmetry of complex-valued random vectors}

\author{Norbert Henze\\
Institute of Stochastics, \\
Karlsruhe Institute of Technology (KIT) \\
Englerstr. 2, D-76133 Karlsruhe \\
\texttt{Norbert.Henze@kit.edu}\\
\And  Pierre Lafaye de Micheaux\\
School of Mathematics and Statistics, University of New South Wales, Sydney, Australia\\
Institut Desbrest d’Epidémiologie et de Santé Publique, IDESP UMR UA11 INSERM - Univ. Montpellier, France\\
AMIS, Université Paul Valéry Montpellier 3, France\\
\texttt{lafaye@unsw.edu.au}\\
\And  Simos G. Meintanis\\
Department of Economics,\\
National and Kapodistrian University of Athens,\\
Athens, Greece \\ and \\ Pure and Applied Analytics \\ North--West University, Potchefstroom, South Africa\\
\texttt{simosmei@econ.uoa.gr}\\
}

\begin{document}

\date{\today}
\maketitle


\begin{abstract}
We propose tests for the null hypothesis that the law of a complex-valued random vector is circularly symmetric.
The test criteria are formulated as $L^2$-type criteria based on empirical characteristic functions,
 and they are convenient from the computational point of view. Asymptotic as well as Monte-Carlo results are presented. Applications on real data are also reported. An R package called \texttt{CircSymTest} is available from the authors.
\end{abstract}

\section{Introduction}\label{secintro}
Let $ Z = (Z^{(1)},\ldots,Z^{(d)})^\top$ be a $\mathbb{C}^d$-valued random (column) vector, where $d \ge 1$ is a fixed integer,
and $\top$ denotes transposition. Moreover, let $ Z_1,  Z_2, \ldots $ be a sequence of independent   and identically distributed (i.i.d.) copies of $ Z$.
We assume that all random vectors are defined on a common suitable probability space $(\Omega,{\cal A}, \mathbb{P})$.
Writing $\edist$ for equality  in distribution, and putting
 $\Theta := [-\pi,\pi)$, we propose and study a test of the hypothesis that the distribution of $ Z$ is (weakly) {\em circularly symmetric}, i.e.,
\begin{equation}\label{NULL}
H_0:  Z \edist \e^{\ii \vartheta}  Z \quad \text{for each} \ \vartheta \in \Theta,
\end{equation}
against general alternatives, on the basis of $ Z_1,\ldots,  Z_n$.

The notion of circular symmetry, as well as its generalization to (complex) elliptical symmetry,
has numerous applications in engineering, and particularly in signal processing. Detailed discussions may be found in \cite{PI:1994},
\cite{OTK:2012}, and \cite{LA:2017}, Chapter 24). We also refer to these publications for the basic definitions and properties
 of complex-valued random variables and random vectors employed herein; see also \cite{DLM:2016}.  On the other hand,
for the related notions of reflective and spherical symmetry, as well as for other forms of symmetry of random vectors in $\RR^d$
and corresponding goodness-of-fit tests,  the reader is referred to \cite{HKM:2003}, \cite{MGW:2012}, and \cite{HHM:2014}.
\textcolor{black}{Typically, estimation in connection with and testing for circular symmetry is based on the so-called ``circularity quotient'' which,
for a zero-mean complex random variable $Z\in \mathbb C^1$,  is defined by $\mathbb E(Z^2)/\mathbb E(|Z|^2)$, i.e., by  the ratio of the pseudo-variance and the variance.
 Generalized likelihood ratio tests that involve the circularity coefficient were initially designed with the complex Gaussian distribution in mind, but they
 have subsequently been extended to non-Gaussianity, and there are also robust versions that work with distributions having infinite moments.
 The reader is referred to \cite{AP:1984}, \cite{PB:1997}, \cite{O:2008}, \cite{OK:2009}, \cite{Walden:09}, \cite{OEK:2011}, \cite{OKP:2011}, and \cite{KDM:2014}, among others, for existing tests of circularity.}

Our approach towards the construction of a class of tests of $H_0$  is different from that in the aforementioned papers in that we employ the notion of the characteristic function. Specifically let
\begin{equation*} \label{CFn}
\varphi_{ Z}( z)=\mathbb  E\big[\e^{\ii \re( z^{\rm{H}}  Z)}\big], \quad  z\in \mathbb C^d,
\end{equation*}
denote the characteristic function (CF) of $ Z$,
where $ z^{\rm{H}}=(\bar z_1,...,\bar z_d)$ denotes the transpose conjugate of  $z:=(z_1,...,z_d)^\top$.
Since the CF uniquely determines the distribution of $ Z$,  the hypothesis $H_0$ may be restated in the equivalent form
\begin{equation} \label{h0equiv}
\varphi_{ Z} = \varphi_{{ Z}_\vartheta} \quad  \text{for each} \ \vartheta \in \Theta,
\end{equation}
where $Z_\vartheta$ is shorthand for  $\e^{\ii \vartheta}  Z$. Likewise, we write
${ Z}_{\vartheta,j} := \e^{\ii \vartheta}  Z_j$, $j \ge 1$.  Our test will be based on the empirical CFs
\begin{equation}\label{defphin}
\varphi_n( z) := \frac{1}{n}\sum_{j=1}^n \exp\left(\ii \re( z^{\rm{H}}  Z_j)\right), \qquad
\varphi_{n,\vartheta}( z) := \frac{1}{n}\sum_{j=1}^n \exp\left(\ii \re( z^{\rm{H}}  Z_{\vartheta,j})\right),
\end{equation}
of $ Z_1, \ldots,  Z_n$ and ${ Z}_{\vartheta,1}, \ldots, { Z}_{\vartheta,n}$, respectively.
Since, under $H_0$, $\varphi_n( z) -  \varphi_{n,\vartheta}( z)$ converges to $0$ $\mathbb{P}$-almost surely as $n \to \infty$ for each
$ z$ and each $\vartheta$, it makes sense to reject $H_0$ for large values of the weighted $L^2$-statistic
\begin{equation}\label{TSn}
T_{n} = n \int_\Theta \int_{\mathbb{C}^d} \big{|} \varphi_n( z)  -  \varphi_{n,\vartheta}( z) \big{|}^2 \, \gamma( z,\vartheta) \, {\rm d}{ z}{\rm d}\vartheta,
\end{equation}
where $\gamma: \mathbb{C}^d \times \Theta \to \mathbb{R}$ is a suitable non-negative weight function.

Although this approach is appealing from a theoretical point of view, we have to impose restrictions on the weight
function $\gamma$ to make a test of $H_0$ based  on $T_{n}$ feasible in practice. Furthermore, we will see that
we can entirely dispense with complex numbers. To this end, we put $ Z =:  X + \ii  Y$, where $ X = {\rm Re}  Z$
and $ Y = {\rm Im}  Z$, and taking real and imaginary parts is understood to be componentwise, i.e.,
${\rm Re}  z := ({\rm Re} z_1, \ldots, {\rm Re} z_d)^\top$, if $ z = (z_1,\ldots,z_d)^\top \in \mathbb{C}^d$,
and likewise for the imaginary part. In the same way, we define $ Z_j =:  X_j + \ii  Y_j$ and
$ Z_{\vartheta,j} =:  X_{\vartheta,j} + \ii  Y_{\vartheta,j}$, $j \ge 1$. Moreover, \textcolor{black}{putting}
\[
\textcolor{black}{\ri(z) := ( \re z_1,\ldots,\re z_d, \im z_1,\ldots,\im z_d)^\top, \qquad z = (z_1,\ldots.z_d)^\top \in \mathbb{C}^d, }
\]
\textcolor{black}{let $W:= \ri(Z)$, $W_j := \ri(Z_j)$, $W_{\vartheta,j} := \ri(Z_{\vartheta,j})$ and $W_\vartheta :=  \ri(Z_\vartheta)$,}
 where $ Z_\vartheta =:  X_\vartheta + \ii  Y_\vartheta$.
Notice that $ W,  W_\vartheta,  W_j$ and $ W_{\vartheta,j}$ take values in $\mathbb{R}^{2d}$. Now, \textcolor{black}{letting}
\begin{equation}\label{defsz}
 s := ({\rm Re}z_1,\ldots,{\rm Re}z_d,{\rm Im}z_1,\ldots,{\rm Im}z_d)^\top
\end{equation}
\textcolor{black}{for the sake of brevity},
then $\re( z^{\rm{H}}  Z) =  s^\top  W$, and straightforward calculations (using $|z|^2 = \bar z z$ for $z \in \mathbb{C}$) yield
\begin{equation}\label{phinmq}
 \big{|} \varphi_n( z)\!  -\!  \varphi_{n,\vartheta}( z) \big{|}^2
 =  \frac{1}{n^2} \sum_{j,k=1}^n \! \left(\cos( s^\top  W_{j,k}) \! - \! 2 \cos( s^\top  W_{j,\vartheta k})
\! + \! \cos(s^\top  W_{\vartheta j,\vartheta k})    \right),
\end{equation}
where
\begin{equation}\label{defofws}
 W_{j,k} :=  W_j\! - W_k, \quad  W_{j,\vartheta k} :=  W_j\! -\!  W_{\vartheta,k},   \quad W_{\vartheta j,\vartheta k} :=  W_{\vartheta,j}\! -\!  W_{\vartheta,k}.
 \end{equation}

In terms of $W$ and $W_\vartheta$, the transformation \textcolor{black}{$ Z \mapsto \exp(\ii \vartheta)  Z$} is equivalent to
$ W \mapsto  W_\vartheta := M_\vartheta    W$, where -- denoting by ${\rm I}_d$ the unit matrix of order $d$ --
$M_\vartheta$ is the ($2d\times 2d$)-matrix
\begin{equation}\label{dmatrix}
 M_\vartheta=\begin{pmatrix} \cos\vartheta \:  {\rm{I}}_d & -\sin\vartheta \: {\rm{I}}_d \\
 \sin\vartheta \: {\rm{I}}_d   &  \cos\vartheta \:  {\rm{I}}_d   \end{pmatrix}.
\end{equation}
If $\varphi_{ U}^\ast (t) := \mathbb{E}(\exp(\ii t^\top  U))$, $t \in \mathbb{R}^{2d}$, denotes the
CF of a $\mathbb{R}^{2d}$-valued random vector $ U$, then  \eqref{h0equiv} holds if and only if
$\varphi_{ W}^\ast = \varphi_{ W_\vartheta}^\ast$  for each $\vartheta \in \Theta$.
Notice that
\[
\varphi_n^\ast( s) := \frac{1}{n} \sum_{j=1}^n \exp\left(\ii  s^\top  W_j\right), \qquad \varphi_{n,\vartheta}^\ast( s) := \frac{1}{n} \sum_{j=1}^n \exp\left(\ii  s^\top  W_{\vartheta,j}\right)
\]
are the 'real world counterparts' of $\varphi_n( z)$ and $\varphi_{n,\vartheta}( z)$ figuring in \eqref{defphin}, and that, in view of \eqref{defsz},
$|\varphi_n^\ast( s) -  \varphi_{n,\vartheta}^\ast( s)|^2$ coincides with the right hand side of \eqref{phinmq}.

Regarding the feasibility of a test of $H_0$ based on $T_{n}$, we assume that $\gamma$ figuring in \eqref{TSn} has the form
\[
\gamma( z,\vartheta) \ = \ w\left(\| z\|_\mathbb{C}\right), \qquad  z \in \mathbb{C}^d, \ \vartheta \in \Theta,
\]
where $w:[0,\infty) \to (0,\infty)$ is a measurable, integrable function, and $\| z\|^2_\mathbb{C} = \sum_{j=1}^d |z_j|^2$.
Thus, the integration with respect to $\vartheta$ is in fact integration with respect to the uniform distribution, since
the factor $(2\pi)^{-1}$ is unimportant. We will elaborate on the function $w$ in the next section.
If $\|\cdot \|$ denotes the Euclidean norm on $\mathbb{R}^{2d}$, then (recall \eqref{defsz}!)
we have $\| z\|_\mathbb{C} = \| s\|$, and the test statistic $T_{n}$ takes the form
\begin{equation}\label{TSn2}
T_{n} = n \int_\Theta \int_{\mathbb{R}^{2d}} |\varphi_n^\ast( s) -  \varphi_{n,\vartheta}^\ast( s)|^2 w\left(\| s\|\right) {\rm d} s{\rm d} \vartheta.
\end{equation}
Although $T_n$ depends on the weight function $w(\cdot)$, this dependence will only be made explicit (i.e., we write $T_{n,w}$) if necessary.

We stress that, at least in principle, also a statistic analogous to that of the Kolmogorov--Smirnov test, i.e.,
\[
\tau_{n,d} : = \sqrt{n} \sup_{\vartheta \in [-\pi,\pi)} \sup_{s\in \mathbb{R}^{2d}} |\varphi_n^\ast( s) -  \varphi_{n,\vartheta}^\ast( s)|
\]
would be an option.  However,  while $\tau_{n,d}$ is by all means a reasonable statistic,
it does not enjoy the full computational reducibility that $T_{n}$ exhibits, since the supremum
 figuring above must be computed as a maximum after discretization.
 While such an approach may be feasible for low dimensions, it certainly becomes problematic as $d$ increases.

\begin{remark} \label{rem1} In our setting, there is the more general notion of {\em unitary invariance},
 which requires $ Z \edist  C  Z$ for each matrix $ C\in \mathbb C^{d\times d}$
 such that $ C^{\rm{H}}= C^{-1}$, where $ C^{\rm{H}}$ denotes the transpose conjugate of $ C$.
The class of distributions that are unitary invariant in $\mathbb C^d$
coincides with the class of complex spherical distributions
in the same dimension, which in turn is equivalent to spherical symmetry of the $\RR^{2d}$-valued joint
vector $(X^\top Y^\top)^\top$ of real and imaginary parts of $ Z$.
Notice that unitary invariance and circular symmetry coincide in the special case $d=1$, since the class of orthogonal  $(2\times2)$-matrices is exhausted by two types of matrices:
One type forms the subclass  ${ \cal{M}}_\vartheta:=\{M_\vartheta, \vartheta\in\Theta\}$,
with $M_\vartheta$ defined in  \eqref{dmatrix} with $\textrm{I}_1:=1$, while the other type forms the subclass    $\widetilde{\cal{M}}_\vartheta:=\{\widetilde M_\vartheta, \vartheta\in\Theta\}$, where
\[\widetilde M_\vartheta:=
\begin{pmatrix} \cos\vartheta \:   & \sin\vartheta \:  \\
 \sin\vartheta \:   &  - \cos\vartheta \:    \end{pmatrix}.
\]
But clearly $\mathbb R^2 \ni (X \: Y):=W \edist \widetilde W:=(X -Y)^\top$ under spherical symmetry,
and thus $W \edist M_\vartheta W$ implies $W \edist \widetilde M_\vartheta W$, since $\widetilde M_\vartheta W = M_\vartheta \widetilde W$.
 \end{remark}

\begin{remark}
We have already seen in Remark \ref{rem1} that circular symmetry is related to spherical symmetry. Actually, however,
circular symmetry is much closer to the weaker notion of {\em reflective symmetry} of a real-valued random variable. In fact, if $d=1$ then
  we can rephrase (\ref{NULL}) to $aZ\stackrel{{\cal{D}}}{=}Z$  for each $a \in \mathbb{C}$ such that $|a|=1$. Now, if $Z$ and $a$ are real-valued,
   then $a=\pm1$, which reduces to the classical definition of symmetry around the origin, i.e., to  $Z\stackrel{{\cal{D}}}{=}-Z$,
   while if $Z$ and $a$ are complex-valued, then $a=\e^{i\vartheta}$ for some $\vartheta$, which leads to (\ref{NULL}).
We will scrutinize this connection a little further. Recall that if $Z$ is real-valued,  then the CF   $\varphi_Z(t):=\mathbb E(\e^{\ii t Z})$, $t \in \RR$,
  is just the centre of mass of the distribution of $tZ$, after having   wrapped this distribution around the unit circle,
  see  \cite{EPP:1993}. Hence, under a location shift,  $Z_\vartheta=Z+\vartheta$, 
  the distribution of $\e^{\ii zZ_\vartheta}$
  is that of $\e^{\ii zZ}$ rotated by the angle $\vartheta z, \ \vartheta \in \mathbb R$.
  Consequently, as $z$ varies (we are then on a cylinder), the distance of the centre of mass remains fixed, i.e., we have
  $|\varphi_{Z_\vartheta}(z)|=|\varphi_Z(z)|$ for each $z\in\mathbb R$, while due to rotation we have
   ${\rm{Arg}}(\varphi_{ Z_\vartheta}(z))= {\rm{Arg}}(\varphi_{Z}(z))+\vartheta z$, where ${\rm{Arg}}(\cdot)$
  stands for the principal argument of a complex number. On the other hand,  if $Z$ is complex-valued, i.e.,  when we are already
  in the complex plane, then $Z_\vartheta=\e^{\ii \vartheta}Z$ is just a rotation of $Z$ itself.
 (Note incidentally that the matrix $M_\vartheta$ defined in (\ref{dmatrix}) is a rotation matrix).
 Consequently,  we have $|Z_\vartheta|=|Z|$ for each $\vartheta$, while due to rotation ${\rm{Arg}}(Z_\vartheta)={\rm{Arg}}(Z)+\vartheta$.
 The above reasoning shows that, since location shifts in $\mathbb R$ map to rotations in $\mathbb C$, it is only natural
 to express these location shifts in the real domain in a  convenient way via identities for the CFs of the corresponding
  random variables $Z$ and $Z_\vartheta$ involved. On the other hand, in the complex domain, the same identities should involve
  the random variables $Z$ and $Z_\vartheta$ themselves, rather than their CFs. Now, assume that $Z$ is real-valued with a distribution
   symmetric around zero. Then the CF of $Z$, as the centre of mass of the wrapped-around-distribution, clearly lies
   on the real axis, i.e., we have  $\varphi_Z(t)={\overline{\varphi}}_Z(t)$ for each $t \in \mathbb R$,
  or equivalently $Z+\vartheta \stackrel{{\cal{D}}}{=}-Z+\vartheta$ for each location shift $\vartheta$, while for  circular symmetry to hold true, the original variable $Z$ and the rotated variable $Z_\vartheta$  must have the same distribution for each rotation $\vartheta$.

\end{remark}
\begin{remark} \label{rem3}
Note that if $\Phi$ is uniform over $[-\pi,\pi)$ and $Z$ is an arbitrary complex-valued random vector independent of $\Phi$,
then $\e^{\ii \Phi}Z$ is circularly symmetric (see \cite{LA:2017}, \S24.3),
a property which \textcolor{black}{also holds if $Z =z$ is fixed, and which} is  analogous to the property that for an arbitrary real-valued random vector $Z$,
the random-signed vector $Z^{(\pm)}:=\pm Z$, with probability 1/2, is symmetrically distributed around zero.
 The last property has been used  for resampling test criteria in the context of testing symmetry of real-valued vectors;
  see \cite{HKM:2003} and \cite{ZHU:2005}, chapter 3. In what follows, the corresponding property for complex-valued random
  vectors  will be the basis for approximating the limit null distribution of the proposed test statistic by the bootstrap.

\end{remark}

The remainder of this work unfolds as follows. In Section~\ref{seccomputation} we deal with computational issues in order to
make the test feasible in practice, while in Section ~\ref{sweight} we investigate the role of the weight function involved in the test statistic $T_n$. Section~\ref{secasymptotics} is devoted to the large-sample behavior of $T_n$.
  In  Section~\ref{SEC5}, we present the results of a simulation study, which has been conducted
   to assess the  finite-sample properties of a resampling version of the new test for circular symmetry. Section~\ref{SEC6} exhibits a
   real-data application. 

%
%
%

\section{Computation of the test statistic}\label{seccomputation}
This section is devoted to computational aspects regarding the test statistic
\be \label{easy}
T_{n} = n \int_\Theta \int_{\mathbb{R}^{2d}} |\varphi_n^\ast( s) -  \varphi_{n,\vartheta}^\ast( s)|^2\,  w\left(\| s\|\right)\,  {\rm d} s{\rm d} \vartheta.
\ee
\textcolor{black} {To this end, a convenient starting point for the choice of the weight function is to consider, apart from a factor, the class of densities
of spherical distributions on $\mathbb R^{2d}$, and to use the fact that these densities -- as well as the corresponding CFs -- depend entirely on the Euclidean norm of
 their argument (see Theorem 2.1 of \cite{FKN:1990}). Hence, due to symmetry,
\be \label{intsinezero}
\textcolor{black}{\int_{\mathbb R^{2d}}  \sin(s^\top u)\,  w(\|s\|)\,  \textrm{d}s = 0,  \qquad u \in \mathbb{R}^{2d},}
\ee
and thus -- apart from a factor -- the corresponding CF is given by
\begin{equation}\label{int}
I(v):=\int_{\mathbb R^{2d}}  \cos(s^\top v)\,  w(\|s\|)\,  \textrm{d}s,  \qquad v \in \mathbb{R}^{2d}.
\end{equation}
Moreover, we have
\begin{equation}\label{INT1n}
I(v)= \Psi(\|v\|) , \qquad v \in \mathbb{R}^{2d};
\end{equation}
where $\Psi$ is labelled the ``characteristic kernel'' of the corresponding spherical distribution.}

Consequently, since $|\varphi_n^\ast( s) -  \varphi_{n,\vartheta}^\ast( s)|^2$ is given by the right hand side of \eqref{phinmq}, it follows that
\begin{equation}\label{shorttstat}
T_n  =  \frac{1}{n} \sum_{j,k=1}^n \int_{-\pi}^\pi \Big{\{} I(W_{j,k}) - 2 I(W_{j,\vartheta k}) + I(W_{\vartheta j,\vartheta k})\Big{\}} \,  {\rm d}\vartheta,
\end{equation}
where $W_{j,k}$,   $W_{j,\vartheta k}$ and $W_{\vartheta j,\vartheta k}$ are given in \eqref{defofws}. Thus in view of \eqref{INT1n} and the fact that $M_\vartheta$ is an orthogonal matrix, we have
\begin{eqnarray}  \nonumber
\|W_{\vartheta j,\vartheta_k}\|^2&=& \|W_{\vartheta,j}-W_{\vartheta,k}\|^2\\ \nonumber
&=&\|M_\vartheta W_{j}-M_\vartheta W_{k}\|^2=\|M_\vartheta W_{j,k}\|^2=\textcolor{black}{W^\top_{j,k}} M^\top_\vartheta M_\vartheta W_{j,k}=\|W_{j,k}\|^2.
\end{eqnarray}
Hence, $I(W_{\vartheta j, \vartheta k})=\Psi(\|W_{\vartheta j, \vartheta k}\|)=\Psi(\|W_{j,k}\|)=I(W_{j,k})$,
which entails further simplification in the computation of $T_n$ figuring in \eqref{shorttstat}, since
\begin{equation*}\label{shorttstat1}
T_n  =  \frac{1}{n} \sum_{j,k=1}^n \Big{\{} 4 \pi \Psi(\|W_{j,k}\|) - 2 \int_{-\pi}^\pi \Psi(\|W_{j, \vartheta k}\|) {\rm d}\vartheta \Big{\}}.
\end{equation*}
As for the integral occurring above, notice that
\begin{eqnarray*}
\|\textcolor{black}{W_{j,\vartheta k}}\|^2 & = & \|W_j- M_\vartheta W_k\|^2\\
& = & \|W_j\|^2 + \|W_k\|^2 - 2 W_j^\top M_\vartheta W_k.
\end{eqnarray*}
Moreover, straightforward calculations yield $W_j^\top M_\vartheta W_k = C_{j,k}\cos \vartheta + S_{j,k}\sin \vartheta$, where
\[
C_{j,k} = X_j^\top X_k + Y_j^\top Y_k, \quad S_{j,k} = X_k^\top Y_j -X_j^\top Y_k.
\]
\textcolor{black}{If we choose the spherical Gaussian density as weight function in \eqref{int},  then the resulting kernel is
\begin{equation} \label{gaussian}
\Psi(\xi) = \e^{-\lambda \xi^2}, \quad \xi \ge 0,  \ \lambda>0,
\end{equation}
where $2\lambda$ is the componentwise variance of the corresponding Gaussian vector}. Furthermore, if we use the fact that
\[
\int_0^{2\pi} \e^{p\cos \vartheta + q \sin \vartheta} \, {\rm d} \vartheta = 2\pi I_0\left(\sqrt{p^2+q^2}\right),
\]
where
\[
I_0(t) = \sum_{k=0}^\infty \frac{t^{2k}}{4^k k!^2}, \quad t \in \mathbb{R},
\]
is the modified Bessel function of the first kind of order $0$ (see  \cite{GR:1994}, \S3.93, equation 3.937 2.),
then $T_n$ takes the form
\be \label{test_form}
T_{n,\lambda} & = &\frac{4\pi}{n} \sum_{j,k=1}^n \Big{[}\e^{-\lambda(\|W_j\|^2+ \|W_k\|^2)} \Big{\{}\e^{2\lambda W_j^\top W_k} -  I_0\left(2\lambda\sqrt{C_{j,k}^2 + S_{j,k}^2} \right)  \Big{\}}      \Big{]}, \ee
where we have made the dependence of $T_n$ on $\lambda$ explicit. 

\globalcolor{black}

\section{The role of the weight function}\label{sweight}


The choice of the characteristic kernel $\Psi$ figuring in \eqref{INT1n} is in part motivated by computational convenience, with simple kernels preferred to more complicated ones.
In this connection, the kernel figuring in (\ref{gaussian}), which corresponds to the Gaussian density, can be generalized if we adopt the spherical stable density
with resulting kernel $\Psi(\xi) := \exp(-\lambda \xi^\mu)$, $\lambda>0$, $\mu \in (0,2]$, for which $\mu=2$ is the boundary Gaussian case; see \cite{NO:2013}.
The existence of alternative kernel choices resulting from different values of $\mu\in(0,2]$  renders the test statistic sufficiently flexible with respect to power performance.
In finite samples, some evidence of this flexibility is provided by \cite{CMZ:2019}. Notwithstanding the significance of the choice of a kernel,
the impact of the weight parameter $\lambda>0$ on the performance of the test may be even more important.  In this regard, we note that a ``good choice'' for $\lambda$
leads to a non-trivial analytic problem, the so-called ``eigenvalue problem'', for which explicit solutions are rarely known; see \cite{Tenreiro:09},  \cite{Tenreiro:11}, \cite{LMR:14}, \cite{HR:19} for  recent contributions. This line of research has led to data-dependent choices for $\lambda$, such as those suggested by  \cite{AS:15} and \cite{Tenreiro:19}. Data-dependent choices, however,
 are tailored to the much more restricted context of testing goodness-of-fit for parametric families of distributions,  so we do not make specific use of such methods here.

 In what follows, we provide some insight into the role of the weight parameter $\lambda$ for the Gaussian kernel figuring in \eqref{gaussian}.
  To this end, we consider the inner integral in \eqref{TSn2}. Due to symmetry (see \eqref{intsinezero}), straightforward algebra yields
 \begin{eqnarray} \label{ttheta}
T_{n,\lambda}(\vartheta) &:=&  \int_{\mathbb{R}^{2d}} |\varphi_n^\ast( s) -  \varphi_{n,\vartheta}^\ast( s)|^2  {\rm e}^{-\lambda \| s\|^2}  {\rm d} s \\ \nonumber 
  \! \! & \! \! = \! \! & \! \!   \int_{\mathbb{R}^{2d}} \left(\frac{1}{n} \sum_{j=1}^n \left\{ \cos s^\top W_j+\sin s^\top W_j-\cos s^\top  W_{\vartheta,j}-\sin s^\top W_{\vartheta,j}\right\}\right)^2  {\rm e}^{-\lambda \| s\|^2}  {\rm d} s
  \\  \nonumber
 \! \! & \! \! = \! \! & \! \!  \int_{\mathbb{R}^{2d}} \! \left(\! \frac{1}{n} \sum_{j=1}^n \! \left \{ \frac{s^\top \! W_j\! -\! s^\top \! W_{\vartheta,j}}{1!}-
 \frac{(s^\top \! W_j)^2\! -\! (s^\top \! W_{\vartheta,j})^2}{2!}-
  \frac{(s^\top \! W_j)^3\! -\! (s^\top \! W_{\vartheta,j})^3}{3!} + \cdots   \right\}\! \! \right)^2  \! {\rm e}^{-\lambda \| s\|^2}  {\rm d} s
  \\
 \! \! & \! \! = \! \! & \! \!    \label{exp}
 \int_{\mathbb{R}^{2d}} \left(\frac{1}{n} \sum_{j=1}^n \frac{s^\top( W_j-W_{\vartheta,j})}{1!}\right)^2 {\rm e}^{-\lambda \| s\|^2}  {\rm d} s \\ \nonumber
 &-& 2   \int_{\mathbb{R}^{2d}}
 \left(\frac{1}{n} \sum_{j=1}^n \frac{s^\top( W_j-W_{\vartheta,j})}{1!}\right) \left( \frac{1}{n} \sum_{j=1}^n \frac{(s^\top W_j)^2-(s^\top W_{\vartheta,j})^2}{2!}\right) {\rm e}^{-\lambda \| s\|^2}  {\rm d} s
\\ \nonumber
&+& \int_{\mathbb{R}^{2d}} \left(\frac{1}{n} \sum_{j=1}^n \frac{(s^\top W_j)^2-(s^\top W_{\vartheta,j})^2}{2!}\right)^2 {\rm e}^{-\lambda \| s\|^2}  {\rm d} s \\ \nonumber
&-& 2   \int_{\mathbb{R}^{2d}}
 \left(\frac{1}{n} \sum_{j=1}^n \frac{s^\top( W_j-W_{\vartheta,j})}{1!}\right) \left( \frac{1}{n} \sum_{j=1}^n \frac{(s^\top W_j)^3-(s^\top W_{\vartheta,j})^3}{3!}\right) {\rm e}^{-\lambda \| s\|^2}  {\rm d} s+ \cdots
   \end{eqnarray}
Here, the last two equations result from the approximations $\sin x=x-(x^3/3!)+\cdots$  and  $\cos x=1-(x^2/2!)+\cdots \ $, and by expansion of the square and grouping 
according to increasing powers. It is already transparent from the expansion in  \eqref{exp} that the test statistic incorporates  empirical moments of contrasts between 
``projections'' of the observations $(W_j, \ j\geq 1)$ on the one hand and corresponding  projections of $(W_{\vartheta,j}, \ j\geq 1)$ on the other hand, as well as empirical moments
 of powers of these contrasts. This fact is quite intuitive, since $W\edist W_\vartheta$ under the null hypothesis of circular symmetry, and consequently these empirical moments should be close to zero.

In order to further scrutinise the role of the weight function, we carry out the multiplications of the sums indicated in \eqref{exp} and integrate term by term, which gives
 \begin{eqnarray} \label{tay}
 T_{n,\lambda}(\vartheta) &=&   \frac{1}{n^2} \sum_{j,k=1}^n \int_{\mathbb{R}^{2d}} \frac{s^\top( W_j-W_{\vartheta,j})}{1!}\frac{s^\top( W_k-W_{\vartheta,k})}{1!} {\rm e}^{-\lambda \| s\|^2}  {\rm d} s
 \\ \nonumber
 &-&   \frac{2}{n^2}   \sum_{j,k=1}^n  \int_{\mathbb{R}^{2d}} \frac{s^\top( W_j-W_{\vartheta,j})}{1!}  \frac{(s^\top W_k)^2-(s^\top W_{\vartheta,k})^2}{2!} {\rm e}^{-\lambda \| s\|^2}  {\rm d} s
 \\ \nonumber
 &-&   \frac{2}{n^2}   \sum_{j,k=1}^n  \int_{\mathbb{R}^{2d}} \frac{s^\top( W_j-W_{\vartheta,j})}{1!}
 \frac{(s^\top W_k)^3-(s^\top W_{\vartheta,k})^3}{3!} {\rm e}^{-\lambda \| s\|^2}  {\rm d} s \\   \nonumber
  \\   \nonumber
&+&  \frac{1}{n^2} \sum_{j,k=1}^n \int_{\mathbb{R}^{2d}} \frac{(s^\top W_j)^2-(s^\top W_{\vartheta,j})^2}{2!}\frac{(s^\top W_k)^2-(s^\top W_{\vartheta,k})^2}{2!} {\rm e}^{-\lambda \| s\|^2}  {\rm d} s+\cdots
\\ \nonumber 
&=:&\frac{1}{n^2} \sum_{j,k=1}^n\left( K^{(2)}_{j,k}-2  \Lambda^{(3)}_{j,k}-2M^{(4)}_{j,k}+N^{(4)}_{j,k}+\cdots\right)
 \end{eqnarray}
 (say). Here, the superscripts denote the total of exponents involved in the projections occurring  in each integral. 
 The expansion in \eqref{tay} involves weighted type  $V$-statistics incorporating the aforementioned contrasts, and it underlines the fact  that the weight function serves the purpose of 
 assigning a specific functional form on the weights of the contrasts. Arguing
analogously to \eqref{intsinezero}, it follows that
\[
\Lambda^{(3)}_{j,k}=0,
\]
 since the integrand is an odd function, and the same holds for any subsequent integral with an odd-numbered superscript.   Furthermore, after some extra algebra one may show that
\[
K^{(2)}_{j,k}={\cal{O}}(\lambda^{-d-1}),  \ \  M^{(4)}_{j,k}, N^{(4)}_{j,k}=   {\cal{O}}(\lambda^{-d-2}).
\]
More generally,  for each  subsequent non-vanishing integral $\Xi$ (say), we have
\[
\Xi^{(2p)}_{j,k}={\cal{O}}(\lambda^{-d-p}),  \ \lambda\to\infty, \ p=3, 4,\ldots 
\]
Thus, the role of the weight parameter $\lambda$ is clearer now: It regulates the rate of decay  by which each power of contrasts enters into the test statistic. 
Specifically, for a fixed value of $\lambda$, higher powers progressively receive less weight, and with increasing $\lambda$ the effect of these powers on the test statistic diminishes.
 This fact is quite intuitive if we take into account that such high powers are more prone to statistical error. 
 This kind of behaviour apparently calls for a compromise between values of $\lambda$ that are too close to the origin and lead to statistical (as well as numerical) instability, and large values that lead to a test based only on a  few contrasts of low order. In this connection and for fixed sample size $n$, the ultimate effect of taking $\lambda$ large is expressed by the following limit value \[
\frac{2}{\pi^{d}}\lim_{\lambda \to \infty}\left[ \lambda^{d+1} T_{n,\lambda}(\vartheta)\right]=2(1-\cos\vartheta)\|\overline W\|^2,
\]
which shows that such extreme weighting leads to a test criterion that rejects the null hypothesis of circular symmetry for a large value of $\|\overline W\|$, where $\overline W$ stands for the
 componentwise sample mean of the observation vectors $(W_j, \ j=1,...,n)$. Clearly, this value should be close to zero, and it should asymptotically vanish for large sample size $n$, under circular symmetry.

\globalcolor{black}

\section{Asymptotic results}\label{secasymptotics}
In this section, we derive asymptotic results for $T_{n}$ defined in \eqref{TSn2}, both under the hypothesis $H_0$ as well as under alternatives.
To this end, we let
\[
{\rm CS}^{+}(\xi) := \cos(\xi) + \sin(\xi), \qquad \xi \in \mathbb{R},
\]
and put
\[
\ell( w,  s,\vartheta) := {\rm CS}^{+}( s^\top  w) - {\rm CS}^{+}( s^\top  w_\vartheta), \qquad  w,  s \in \RR^{2d}, \ \vartheta \in \Theta,
\]
where $ w_\vartheta = M_\vartheta  w$, and $M_\vartheta$ is defined in \eqref{dmatrix}. Furthermore, we write
\begin{equation}\label{defvn}
V_n( s,\vartheta) = \frac{1}{\sqrt{n}} \sum_{j=1}^n \ell( W_j, s,\vartheta).
\end{equation}
Using \eqref{intsinezero}, addition theorems for the
sine and the cosine functions, and some calculations show that
\[
T_{n} = \int_\Theta \int_{\RR^{2d}} V^2_n( s,\vartheta)  w(\| s\|) \, {\rm d} s {\rm d}\vartheta.
\]
Let $\mathbb{H} = L^2(\RR^{2d}\times \Theta,{\cal B}(\RR^{2d}\times \Theta),w(\| s\|){\rm d} s{\rm d}\vartheta)$ denote the separable Hilbert space of
(equivalence classes of) measurable functions $f:\RR^{2d}\times \Theta \to \RR$ that are square integrable with respect to
$w(\|  s\|) {\rm d} s{\rm d}\vartheta $. The inner product and the norm in $\mathbb{H}$ will be denoted by
\[
\langle f,g \rangle_\mathbb{H} = \int_\Theta \int_{\RR^{2d}} f( s,\vartheta)g( s,\vartheta)  w(\| s\|) \, {\rm d} s {\rm d}\vartheta, \qquad \|f\|_\mathbb{H} = \langle f,f \rangle_\mathbb{H}^{1/2},
\]
respectively. Notice that
$\ell( W_j,\cdot ,\cdot), \ j \ge 1$, is an i.i.d. sequence of random elements of $\mathbb{H}$ satisfying
\begin{equation}\label{secmom}
\mathbb{E} \|\ell( W_1,\cdot,\cdot)\|_\mathbb{H}^2 = \int_\Theta \int_{\RR^{2d}} \mathbb{E} \big{[} \ell^2( W_1( s,\vartheta))\big{]}  w(\| s\|) \, {\rm d} s {\rm d}\vartheta < \infty.
\end{equation}
Our first result is an almost  sure limit for $n^{-1}T_{n}$.

\begin{theorem}\label{theo1} Without any restriction on the distribution of $ Z$,  the statistic $T_{n}$ in \eqref{TSn2} satisfies
\[
\lim_{n\to \infty} \frac{T_{n}}{n}  =  \Delta \quad \mathbb{P}\text{-almost surely},
\]
where
\begin{equation}\label{defvondelta}
\Delta = \int_\Theta \int_{\mathbb{R}^{2d}} |\varphi^\ast_{ W}( s) - \varphi^\ast_{ W_{\vartheta}}( s)|^2 w(\| s\|) {\rm d} s {\rm d}\vartheta.
\end{equation}
\end{theorem}

\noindent {\sc Proof.} By the strong law of large numbers in Hilbert spaces \textcolor{black}{(see, e.g., \cite{BO:2000}, Theorem 2.4)},
$n^{-1} \sum_{j=1}^n \ell( W_j,\cdot,\cdot) \rightarrow \mathbb{E}\big{[}\ell( W,\cdot,\cdot)\big{]}$
$\mathbb{P}$-a.s. and thus
\[
\frac{T_{n}}{n} \rightarrow \big{\|}\mathbb{E}\big{[}\ell( W,\cdot,\cdot)\big{]}\big{\|}_\mathbb{H}^2 = \int_\Theta \int_{\mathbb{R}^{2d}} \left(\mathbb{E}[\ell( W, s,\vartheta)]\right)^2
w(\| s\|) \, {\rm d} s {\rm d}\vartheta
\]
as $n \to \infty$ $\mathbb{P}$-a.s. Using symmetry arguments, it is readily seen that the last expression equals $\Delta$ figuring in \eqref{defvondelta}.
Notice also that
\[
\Delta = \int_\Theta \int_{\mathbb{C}^{d}} |\varphi_{ Z}( z) - \varphi_{ Z_{\vartheta}}( z)|^2 w(\| z\|) \, {\rm d} z {\rm d}\vartheta. \ \bewend
\]

In view of Theorem \ref{theo1}, the non-negative quantity $\Delta$, which depends on the weight function $w$, defines  the 'distance to symmetry' in the sense of $H_0$ of the underlying distribution
of $ Z$, and we have $\Delta =0$ if and only if $H_0$ holds.

Regarding the asymptotic null distribution of $T_{n}$ as $n \to \infty$, we have the following result.

\begin{theorem}\label{theo2} If $H_0$ holds, there is a centred Gaussian random element $V$ of $\mathbb{H}$ with covariance kernel $K( s,\vartheta; t,\eta) = \mathbb{E}[V( s,\vartheta)V( t,\eta)]$ given by
\begin{equation}\label{defcovfkt}
K( s,\vartheta; t,\eta) = \mathbb{E}\big{[} \ell( W, s,\vartheta) \ell( W, t,\eta)\big{]}, \quad  s, t \in \mathbb{R}^{2d}, \vartheta_1, \eta \in \Theta,
\end{equation}
such that $V_n \verk V$, where $V_n$ is given in \eqref{defvn}.
\end{theorem}

Since $T_{n} = \|V_n\|_\mathbb{H}^2$, the continuous mapping theorem yields the following corollary.
\begin{corollary} Under $H_0$, we have
\[
T_{n} \verk T_{\infty} := \textcolor{black}{\|V\|_{\mathbb{H}}^2 = } \int_\Theta \int_{\RR^{2d}} V^2( s,\vartheta)  w(\| z\|) \, {\rm d} z {\rm d}\vartheta,
\]
where $V$ is the Gaussian random element of Theorem \ref{theo2}.
\end{corollary}

\medskip

\noindent {\sc Proof of Theorem \ref{theo2}}.  Under $H_0$, the summands figuring in \eqref{defvn} are centred random elements
of $\mathbb{H}$ satisfying  \eqref{secmom}. By the Central limit theorem for i.i.d. random elements in Hilbert spaces, see, e.g., Theorem 2.7 in \cite{BO:2000}), 
there is  a centred Gaussian random element $V$ of $\mathbb{H}$ with covariance function $K$ given in \eqref{defcovfkt},
such that $V_n \verk V$ as $n \to \infty$.  \bewend

\medskip
Since both the finite-sample and the limit null distribution of $T_{n}$ depend on the underlying unknown distribution
of $W$, we suggest the following bootstrap procedure to carry out the test in practice.
Independently of $W,W_1,W_2, \ldots$, let $\Phi, \Phi_1,\Phi_2,\ldots $ be a sequence of i.i.d. random random variables
such that the distribution  of $\Phi$ is uniform over $\Theta = [-\pi,\pi)$. We assume that all random variables are defined on a
common probability space $(\Omega,{\cal A},\PP)$. For given $\omega \in \Omega$, the bootstrap procedure
conditions on the realizations $w_1=W_1(\omega), \ldots, w_n = W_n(\omega)$.
The rationale of this procedure is as follows: Given $w_1,\ldots, w_n$, we have to generate a distribution that satisfies $H_0$.
\textcolor{black}{Since the distribution  of $M_{\Phi_j}w_j$, $j=1,\ldots,n,$ is circularly symmetric,} the significance of the observed value $T_n(w_1,\ldots,w_n)$ of the test statistic
 should be jugded with respect to the distribution of $T_{n}(M_{\Phi_1} w_1, \ldots, M_{\Phi_n} w_n)$.
The latter distribution  can be estimated as follows:
   Choose a large number $B$ and, conditionally on
$ W_1=w_1,\ldots,  W_n=w_n$, generate independent copies
\[
T_{n}^{(b)} := T_{n}(M_{\Phi_{b,1}} w_1, \ldots, M_{\Phi_{b,n}} w_n), \qquad b = 1,\ldots, B,
\]
where $\Phi_{b,j}$, $b \in \{1,\ldots,B\}$, $j \in \{1,\ldots,n\}$,  are i.i.d. with a uniform distribution on $[-\pi,\pi)$. The critical value
for a test of level $\alpha$ based on $T_{n}$ is then the upper $(1-\alpha)$-quantile of the empirical distribution of
$T_{n}^{(b)}$, $b=1,\ldots, B$. The following result shows the asymptotic validity of this bootstrap procedure.

\begin{theorem}\label{theo3} Assume that $H_0$ holds.  For $w_1, \ldots, w_n \in \RR^{2d}$, let
\[
V_n^\ast(s,\vartheta) := \frac{1}{\sqrt{n}} \sum_{j=1}^n \ell(M_{\Phi_j}w_j,s,\vartheta), \quad s \in \RR^{2d}, \vartheta \in \Theta,
\]
and put $T_n^\ast := \|V_n^\ast\|_\mathbb{H}^2$. For $\PP$-almost all sample sequences $W_1(\omega)=w_1, W_2(\omega)=w_2, \ldots $, we have
\[
V_n^\ast \verk V \ \text{and} \ T_n^\ast \verk \|V\|_\mathbb{H}^2
\]
as $n \to \infty$, where $V$ is the Gaussian process figuring in the statement of Theorem \ref{theo2}.
\end{theorem}

\noindent {\sc Proof.}
For $w,s,t \in \RR^{2d}$ and $\vartheta, \eta \in \Theta$, let
\begin{equation}\label{deffuncf}
f(w,s,\vartheta,t,\eta) := \BE\big{[} \ell(M_\Phi w,s,\vartheta) \ell(M_\Phi w,t,\eta)\big{]}.
\end{equation}
Let $D \subset \RR^{2d}\times \Theta$ be a countable dense set. From the strong law of large numbers and the fact that a
countable intersection of sets of probability one has probability one, there is a measurable subset $\Omega_0$ of $\Omega$
such that $\PP(\Omega_0)=1$ and
\begin{equation}\label{covconvergence}
\lim_{n\to \infty} \frac{1}{n} \sum_{j=1}^n f(W_j(\omega),s,\vartheta,t,\eta)  = \BE \big{[} \ell(M_\Phi W,s,\vartheta) \ell(M_\Phi W,t,\eta)\big{]}
\end{equation}
for each $\Omega \in \Omega_0$ and each $(s,\vartheta)\in D$ and $(t,\eta) \in D$. Notice that, by the definition of the function $\ell$ and the Lipschitz
continuity of the sine and the cosine function, convergence in \eqref{covconvergence} is in fact for each $(s,\vartheta)$ and $(t,\eta)$ in $\RR^{2d}\times \Theta$.

In what follows, fix $\omega \in \Omega_0$, and let $w_j := W_j(\omega)$, $j \ge 1$. We have
\[
V_n^\ast(s,\vartheta) =  \sum_{j=1}^n V_{n,j}^\ast(s,\vartheta), \quad (s,\vartheta) \in \RR^{2d}\times \Theta,
\]
where $V_{n,j}^\ast(s,\vartheta) :=  n^{-1/2} \ell(M_{\Phi_j}w_j,s,\vartheta)$.
Notice that
\begin{eqnarray*}
\BE \big{[} V_{n,j}^\ast(s,\vartheta) \big{]} & = & \BE[\text{CS}^+(M_{\Phi_j}w_j,s,\vartheta)] - \BE[\text{CS}^+(M_\vartheta M_{\Phi_j} w_j,s,\vartheta)]\\
& = & 0,
\end{eqnarray*}
since $M_{\Phi_j} w_j \edist M_\vartheta M_{\Phi_j}w_j$ ($= M_{\varphi + \Phi_j}w_j$). Thus $V_{n,j}^\ast = V_{n,j}^\ast(\cdot,\cdot)$ is a centred random element of $\mathbb{H}$. Moreover,
we have $\BE  \|V_{n,j}^\ast\|^2_{\mathbb{H}} < \infty$.
To prove that the sequence of random elements $V_n^\ast = V_n^\ast(\cdot,\cdot)$ of $\mathbb{H}$ converges in distribution to the centred Gaussian random element $V$ of $\mathbb{H}$
figuring in Theorem \ref{theo2}, let $\{e_1, e_2, \ldots \}$ be some complete orthonormal subset of $\mathbb{H}$, and let $C_n$ denote the
covariance operator of $V_n^\ast$. According to Lemma 4.2 of \cite{KMM:2000}, we have to show the following:
\begin{enumerate}
\item[(a)] $\lim_{n\to \infty} \langle C_ne_k,e_\ell \rangle_\mathbb{H} = a_{k\ell}$ (say) exists for each $k,\ell \ge 0$.
\item[(b)] $\lim_{n\to \infty} \sum_{k=0}^\infty \langle C_n e_k,e_k \rangle_\mathbb{H} = \sum_{k=1}^\infty a_{kk} < \infty$.
\item[(c)] $\lim_{n\to \infty} L_n(\varepsilon,e_k) =0 $ for each $\varepsilon >0$ and each $k \ge 0$, where\\[1mm]
$L_n(\varepsilon, h) = \sum_{j=1}^n \BE\big{[} \langle V_{n,j}^\ast,h\rangle_\mathbb{H}^2 {\bf 1}\{|\langle V_{n,j}^\ast,h\rangle_\mathbb{H}|> \varepsilon \}\big{]}$, $h \in \mathbb{H}$.
\end{enumerate}
As for (a), let
\[
K_n(s,\vartheta,t,\eta) := \BE\big{[}V_n^\ast(s,\vartheta) V_n^\ast(t,\eta) \big{]}.
\]
Some algebra and symmetry yield
\[
 K_n(s,\vartheta,t,\eta) = \frac{1}{n} \sum_{j=1}^n f(w_j,s,\vartheta,t,\eta),
\]
where $f$ is given in \eqref{deffuncf}. From \eqref{covconvergence} and the fact that
$M_\Phi W \edist W$, we have pointwise convergence $\lim_{n\to \infty} K_n =K$, where $K$ is given in
\eqref{defcovfkt}. Furthermore, putting $\text{D}(s,\vartheta,t,\eta) := w(\|s\|) w(\|t\|) \text{d}s\text{d}\vartheta \text{d}t\text{d}\eta$, dominated convergence  yields
\begin{eqnarray*}
\lim_{n\to \infty} \langle C_ne_k,e_\ell \rangle_\mathbb{H} & = & \lim_{n\to \infty} \iint \iint K_n(s,\vartheta,t,\eta) e_k(s,\vartheta) e_\ell(t,\eta) \, \text{D}(s,\vartheta,t,\eta) \\
& = & \iint \iint  K(s,\vartheta,t,\eta) e_k(s,\vartheta) e_\ell(t,\eta) \, \text{D}(s,\vartheta,t,\eta)\\
& = & \langle Ce_k,e_\ell \rangle_\mathbb{H},
\end{eqnarray*}
where $C$ is the covariance operator of $V$, and
each of the double integrals is over $\RR^{2d}\times  \Theta$. Setting $a_{k \ell} := \langle e_k,e_\ell \rangle_\mathbb{H}$, condition (a) follows. To prove condition (b), notice that,
by monotone convergence, Parseval's inequality and dominated convergence,  we have
\begin{eqnarray*}
\lim_{n\to \infty} \sum_{k=0}^\infty \langle C_n e_k,e_k\rangle_\mathbb{H} & = & \lim_{n\to \infty} \sum_{k=0}^\infty \BE\langle e_k,V_n^\ast \rangle_\mathbb{H}^2\\
& = & \lim_{n\to \infty} \BE \|V_n^\ast\|_\mathbb{H}^2\\
& = & \int_{\RR^{2d}} \int_\Theta \lim_{n\to \infty} K_n(s,\vartheta,s,\vartheta) w(\|s\|) \, \text{d}s \text{d}\vartheta \\
& = & \int_{\RR^{2d}} \int_\Theta  K(s,\vartheta,s,\vartheta) w(\|s\|) \, \text{d}s \text{d}\vartheta \\
& = & \BE \|V\|_\mathbb{H}^2 \\
& = & \sum_{k=0}^\infty a_{kk} < \infty,
\end{eqnarray*}
which shows that condition (b) holds. Finally, observe that
\begin{eqnarray*}
|\langle V_{n,j}^\ast,e_k \rangle_\mathbb{H}| & = & \frac{1}{\sqrt{n}} \bigg{|}\int_{\RR^{2d}}\int_\Theta \ell(M_{\Phi_j}w_j,s,\vartheta) e_k(s,\vartheta) w(\|s\|) \, \text{d}s \text{d}\vartheta \bigg{|}\\
& \le & \frac{1}{\sqrt{n}} \int_{\RR^{2d}} \int_\Theta \big{|} \ell(M_{\Phi_j}w_j,s,\vartheta) e_k(s,\vartheta)\big{|} w(\|s\|) \, \text{d}s \text{d}\vartheta\\
& \le & \frac{1}{\sqrt{n}} \left(\int_{\RR^{2d}}\int_\Theta  \ell^2(M_{\Phi_j}w_j,s,\vartheta) e_k(s,\vartheta) w(\|s\|) \, \text{d}s \text{d}\vartheta \right)^{1/2} \|e_k\|_\mathbb{H}\\
& \le & \frac{4}{\sqrt{n}},
\end{eqnarray*}
since $|\ell| \le 4$. It follows that $\lim_{n\to \infty} L_n(\varepsilon,e_k) =0 $, which entails the validity of (c).
\bewend

We now show that the test statistic $T_{n}$ has an asymptotic normal distribution under fixed alternatives to $H_0$. The
reasoning closely follows \cite{BEH:2017}.

\begin{theorem} Suppose that $H_0$ does not hold. We then have
\[
\sqrt{n}\left(\frac{T_{n}}{n} - \Delta \right) \verk {\rm N}(0,\sigma^2) \quad \text{as} \ n \to \infty,
\]
where
\[
\sigma^2 = 4 \int_{\Theta \times \mathbb{R}^{2d}} \int_{\Theta \times \mathbb{R}^{2d}} \widetilde{K}( s,\vartheta; t,\eta) v( s,\vartheta)v( t,\eta) w(\| s\|)\,
w(\| t\|) {\rm d} s\, {\rm d}  t \, {\rm d}\vartheta \, {\rm d}\eta,
\]
and $\widetilde{K}( s,\vartheta; t,\eta)$ is given in \eqref{kststar}.
\end{theorem}

\noindent {\sc Proof.} Let $\widetilde{V}_n( s,\vartheta) := n^{-1/2}V_n( s,\vartheta)$, where $V_n( s,\vartheta)$ is given in \eqref{defvn},
and put $v( s,\vartheta) := \mathbb{E}[\ell( W, s,\vartheta)]$. Regarded as elements of $\mathbb{H}$, we write
$\widetilde{V}_n$ and $v$. We then have
\begin{eqnarray}\nonumber
\sqrt{n}\left(\frac{T_{n}}{n} - \Delta \right) & = & \sqrt{n}\left(\|\widetilde{V}_n\|_\mathbb{H}^2 - \|v\|_\mathbb{H}^2\right)\\ \nonumber
& = & \sqrt{n} \langle \widetilde{V}_n -v,\widetilde{V}_n + v\rangle_\mathbb{H}\\ \nonumber
& = & \sqrt{n} \langle \widetilde{V}_n -v,2v + \widetilde{V}_n - v\rangle_\mathbb{H}\\ \label{zerl}
& = & 2\langle\sqrt{n}(\widetilde{V}_n -v),v\rangle_\mathbb{H} + \frac{1}{\sqrt{n}} \|\sqrt{n}(\widetilde{V}_n- v)\|_\mathbb{H}^2.
\end{eqnarray}
Now,
\[
\sqrt{n}(\widetilde{V}_n( s,\vartheta) -v( s,\vartheta)) = \frac{1}{\sqrt{n}} \sum_{j=1}^n  \big{\{} \ell( W_j, s,\vartheta) - \mathbb{E}[\ell( W, s,\vartheta)]\big{\}},
\]
and, invoking once more the Central limit  theorem in Hilbert spaces, there is a centred random element $\widetilde{V}$ of $\mathbb{H}$ having covariance kernel
\begin{equation}\label{kststar}
\widetilde{K}( s,\vartheta; t,\eta) := \mathbb{E}\big{[} \ell( W, s,\vartheta) \ell( W, t,\eta)\big{]} - v( s,\vartheta)v( t,\eta), \quad  s, t \in \mathbb{R}^{2d}, \vartheta, \eta \in \Theta,
\end{equation}
such that $\sqrt{n}(\widetilde{V}_n -v) \verk \widetilde{V}$ as $n \to \infty$.
From \eqref{zerl} and Slutski's lemma, it follows that
$n^{-1/2}(T_{n}- \Delta) \verk 2\langle \widetilde{V},v\rangle_\mathbb{H}$. The distribution of $2\langle \widetilde{V},v\rangle_\mathbb{H}$ is the normal distribution ${\rm N}(0,\sigma^2)$. \bewend

\section{Simulations}\label{SEC5}

This section gathers the results of a simulation study regarding the finite sample properties of the new test for circular symmetry. Throughout this section, the number of Monte Carlo replications is set to $M=10,000$, unless indicated otherwise, and the number of bootstrap replications is set to $B=200$. Moreover, the significance level is set to $\alpha=0.05$. We implemented the test in \eqref{test_form} by the equivalent formula
\begin{eqnarray*}
T_{n,\lambda}
& = & \frac{4\pi}{n} \sum_{j,k=1}^n \Big{[}{\rm e}^{-\lambda\|W_j-W_k\|^2}-\frac{1}{\pi}\int_0^{\pi}{\rm e}^{-\lambda(\|W_j\|^2+\|W_k\|^2)+2\lambda\sqrt{C_{j,k}^2+S_{j,k}^2}\cos(t)}\rm{d}t\Big{]},\end{eqnarray*}
which turned out to be numerically more stable.\\

\textcolor{black}{First, in Section~\ref{sec:compleGauss} we consider the simplest case of the scalar zero-mean complex Gaussian random variable with unit variance and circularity quotient  $\rho$. This case enables us to demonstrate empirically that the implementation of our test behaves as intended. }

\textcolor{black}{Next, we investigate various scenarios of non-circular random variables or vectors to  compare the empirical power of our test (for some choices of $\lambda$) to those of some  competitors. Among others, we consider the univariate adjusted Generalized Likelihood Ratio Test (GLRT), the robust GLRT (RobGLRT) and the Wald's type maximum likelihood (WTML) test, as well as the multivariate GLRT (mGLRT); see \cite{O:2008}, \cite{OEK:2011}, \cite{OK:2009}, \cite{OKP:2011}, and \cite{Walden:09}.}

\subsection{The complex Gaussian random variable}\label{sec:compleGauss}

Let $Z=X+\ii Y$ be a zero-mean scalar complex Gaussian random variable with probability density function
\begin{equation}\label{pdf4.5}
p(z) = \frac{1}{\pi C_{zz}\sqrt{1-|\rho|^2}}\exp\left\{-\frac{-|z|^2-\rm{Re}(\rho \bar{z}^2)}{C_{zz}(1-|\rho|^2)}\right\},
\end{equation}
where
$$C_{zz}=\mathbb{E}\{Z\bar{Z}\}=\mathbb{E}[X^2]+\mathbb{E}[Y^2]=\mathbb{V}(X)+\mathbb{V}(Y)
\text{ and }
\mathbb{E}\{ZZ^\top\}=\mathbb{E}\{Z^2\}=\mathbb{V}(X)-\mathbb{V}(Y)+2\ii\E(XY),
$$
and \textcolor{black}{where}
$$
\rho=\rho_x+\ii\rho_y=\frac{\mathbb{V}(X)-\mathbb{V}(Y)}{C_{zz}}+\ii\frac{2\sqrt{\mathbb{V}(X)}\sqrt{\mathbb{V}(Y)}\rho_{xy}}{C_{zz}}
$$
\textcolor{black}{is the circularity quotient defined in \cite{O:2008}.}

We have
$$\mathbb{V}(X)=0.5C_{zz}(1+\rm{Re}(\rho))
\text{ and }
\mathbb{V}(Y)=0.5C_{zz}(1-\rm{Re}(\rho)),$$
and
$$\rho_{xy}=\frac{\rm{Im}\rho}{\sqrt{1-(\rm{Re}\rho)^2}}$$
with $\rho_{xy}=\pm1$ if $\rm{Re}(\rho)=\pm1$.

Note that we can consider the following three cases:
\begin{itemize}
\item $\mathbb{V}(X)=\mathbb{V}(Y)$, and $X$ is uncorrelated with $Y$, in which case $Z$ is proper and circular,
\item $\mathbb{V}(X)\ne \mathbb{V}(Y)$, and $X$ is uncorrelated with $Y$, in which case $Z$ is noncircular,
\item $\mathbb{V}(X)=\mathbb{V}(Y)$, and $X$ is correlated with $Y$, in which case $Z$ is noncircular.
\end{itemize}
We conducted the following simulations. We set $C_{zz}=1$, and we constructed a grid of values of $\rho=r_k\e^{\ii\theta_{\ell}}\in\mathbb{C}$ on the unit disc, with $r_k=k/9$, $k=0,\ldots,9$, and $\theta_{\ell}=\ell(2\pi)/35$, $\ell=0,\ldots,35$. For each value of $\rho$ on this grid, we generated  $M=1,000$ samples of size $n=10,20$ and $50$ from the corresponding probability distribution function in~\eqref{pdf4.5} and computed the empirical power for our test statistic $T_{n,\lambda}$ with $\lambda=1.0$. Figure~\ref{fig:5.5.1} illustrates the results. The empirical power is clearly increasing with the sample size $n$, and it is isotropic, apart from the two locations $(\rho_x,\rho_y)=(\pm1,0)$, which are associated with circularity.

\begin{figure}[H]
    \centering
   \includegraphics[width=5.5cm]{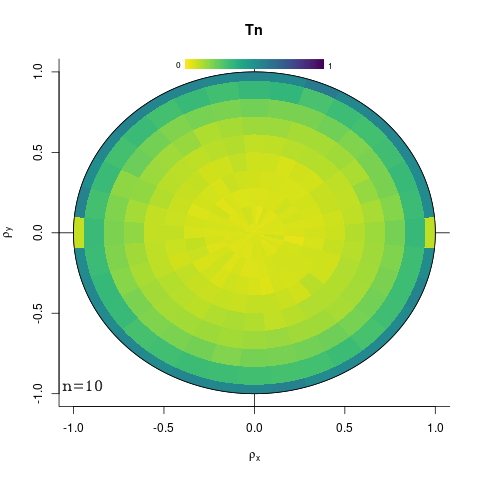}\includegraphics[width=5.5cm]{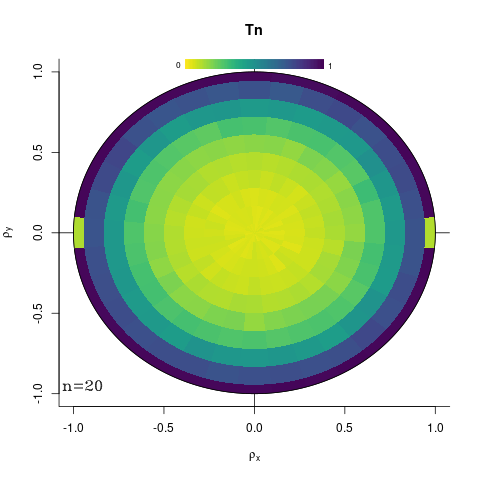}\includegraphics[width=5.5cm]{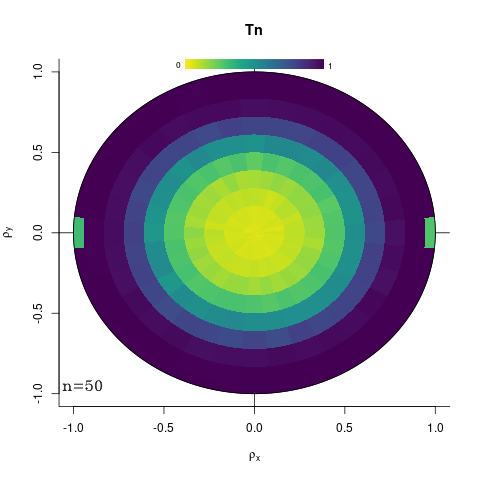}
    \caption{From left to right: empirical power for the test statistic $T_{n,\lambda}$, with $\lambda=1$, for the sample sizes $n=10,20$ and $50$, and the distribution considered in \eqref{pdf4.5}.}
    \label{fig:5.5.1}
\end{figure}

\subsection{Complex Gaussian random vector with identity covariance and varying location}

For a complex Gaussian random vector, circularity is equivalent to the vector being {\it{proper}}. We recall that a complex random vector $Z$ is, by definition,  proper, if the following three conditions are satisfied:
\begin{itemize}
    \item $\E(Z)=0$,
    \item $\mathbb{V}\left(Z^{(1)}\right)<\infty,\ldots,\mathbb{V}\left(Z^{(d)}\right)<\infty$,
    \item $\E[ZZ^\top]=0$.
\end{itemize}

Here, we depart from circularity by allowing the mean of a bivariate complex Gaussian random vector $Z$ (hence $d=2$) to take values increasingly departing from 0, while keeping the two last conditions satisfied. So, we generate random samples of sizes $n=20,50$ and $100$ from a $\textrm{CN}_2((u,u),\rm{I}_2)$ distribution, where $u=0.05k$, $k=0,\ldots,10$. The number of Monte Carlo replications is set to $M=10,000$.

The empirical power results for our test statistic $T_{n,\lambda}$ are presented in Tables~\ref{tab:5.1.1}, \ref{tab:5.1.2} and  \ref{tab:5.1.3}, for $\lambda=0.01$, $\lambda=0.1$ and $\lambda=1.0$, respectively.  \textcolor{black}{Likewise, Table~\ref{tab:5.1.4} shows the corresponding values for the mGLRT of circularity of \cite{OK:2009}. Analogous results for the tests 
$T_1$ and $T_2$ of \cite{Walden:09} that involve the sample correlations between a given complex random vector and its conjugate  
are given in Tables~~\ref{tab:5.1.5} and~\ref{tab:5.1.6}, respectively.} We observe that, \textcolor{black}{when $u=0.00$}, the empirical level is close to the nominal one for all sample sizes considered. Also, as expected, the empirical power increases with the sample size $n$ and with the value of $u=\E(Z^{(1)})=\E(Z^{(2)})$. \textcolor{black}{For these alternatives, it is better to choose a small value of $\lambda$ (i.e., $0.01$) for our test statistic $T_{n,\lambda}$, though all three choices lead to a test that is clearly superior to the mGLRT, $T_1$ and $T_2$ tests.} 

\begin{table}[H]
    \centering
    \caption{Empirical power of $T_{n,\lambda}$ against the alternative  $\textrm{CN}_2((u,u),\rm{I}_2)$  ($\lambda=0.01$, $n=20,50$ and $100$). }
    \label{tab:5.1.1}
    \begin{tabular}{r|rrrrrrrrrrr}
     $u$   &  0.00 & 0.05 & 0.10 & 0.15 & 0.20 & 0.25 & 0.30 & 0.35 & 0.40 & 0.45 & 0.50 \\
        \hline
20  & 0.0542 & 0.0652 & 0.0982 & 0.1659 & 0.2585 & 0.3812 & 0.5321 & 0.6885 & 0.8080 & 0.8942 & 0.9521 \\
50  & 0.0517 & 0.0858 & 0.1753 & 0.3525 & 0.5895 & 0.8202 & 0.9403 & 0.9858 & 0.9985 & 1.0000 & 1.0000 \\
100 & 0.0543 & 0.1075 &  0.3314 &  0.6723 &  0.9097 & 0.9912 & 0.9996 & 1.0000 & 1.0000 & 1.0000 & 1.0000  \\
    \end{tabular}
\end{table}

\begin{table}[H]
    \centering
    \caption{Empirical power of $T_{n,\lambda}$ against the alternative  $\textrm{CN}_2((u,u),\rm{I}_2)$ ($\lambda=0.1$, $n=20,50$ and $100$). }
    \label{tab:5.1.2}
    \begin{tabular}{r|rrrrrrrrrrr}
     $u$   &  0.00 & 0.05 & 0.10 & 0.15 & 0.20 & 0.25 & 0.30 & 0.35 & 0.40 & 0.45 & 0.50 \\
        \hline
20 & 0.0555 & 0.0663 & 0.0992 & 0.1613 & 0.2537 & 0.3746 & 0.5230 & 0.6789 & 0.8009 & 0.8869 & 0.9481  \\
50 & 0.0511 & 0.0845 & 0.1725 & 0.3466 & 0.5797 & 0.8128 & 0.9339 & 0.9838 & 0.9983 & 0.9999 & 1.0000  \\
100 & 0.0526 & 0.1082 & 0.3247 & 0.6623 & 0.9014 & 0.9873 & 0.9995 & 1.0000 & 1.0000 & 1.0000 & 1.0000 \\
    \end{tabular}
\end{table}


\begin{table}[H]
    \centering
        \caption{Empirical power of $T_{n,\lambda}$ against the alternative  $\textrm{CN}_2((u,u),\rm{I}_2)$ ($\lambda=1$, $n=20,50$ and $100$). }
    \label{tab:5.1.3}
    \begin{tabular}{r|rrrrrrrrrrr}
     $u$   &  0.00 & 0.05 & 0.10 & 0.15 & 0.20 & 0.25 & 0.30 & 0.35 & 0.40 & 0.45 & 0.50 \\
        \hline
20 & 0.0529 & 0.0602 & 0.0784 & 0.1160 & 0.1698 & 0.2414 & 0.3402 & 0.4602 & 0.5735 & 0.6849 & 0.7818 \\
50 & 0.0552 & 0.0694 & 0.1226 & 0.2252 & 0.3795 & 0.5909 & 0.7613 & 0.8863 & 0.9575 & 0.9887 & 0.9982 \\
100 & 0.0541 & 0.0871 & 0.2116 & 0.4413 & 0.6973 & 0.9059 & 0.9790 & 0.9978 & 1.0000 & 1.0000 & 1.0000 \\
    \end{tabular}
\end{table}

\globalcolor{black}

\begin{table}[H]
    \centering
        \caption{Empirical power of mGLRT against the alternative  $\textrm{CN}_2((u,u),\rm{I}_2)$ ($n=20,50$ and $100$). }
    \label{tab:5.1.4}
    \begin{tabular}{r|rrrrrrrrrrr}
     $u$   &  0.00 & 0.05 & 0.10 & 0.15 & 0.20 & 0.25 & 0.30 & 0.35 & 0.40 & 0.45 & 0.50 \\
        \hline
20 & 0.0852 & 0.0929 & 0.0945 & 0.0876 & 0.0872 & 0.1009 & 0.1086 & 0.1292 & 0.1562 &  0.1821 & 0.2229 \\
50 & 0.0629 & 0.0629 & 0.0728 & 0.0677 & 0.0731 & 0.0902 & 0.1240 & 0.1610 & 0.2348 &  0.3238 & 0.4454 \\
100 & 0.0531 & 0.0580 & 0.0590 & 0.0642 & 0.0821 & 0.1198 & 0.1797 & 0.2782 & 0.4320 &  0.5996 & 0.7678 \\
    \end{tabular}
\end{table}

\begin{table}[H]
    \centering
        \caption{Empirical power of $T_1$ against the alternative  $\textrm{CN}_2((u,u),\rm{I}_2)$ ($n=20,50$ and $100$). }
    \label{tab:5.1.5}
    \begin{tabular}{r|rrrrrrrrrrr}
     $u$   &  0.00 & 0.05 & 0.10 & 0.15 & 0.20 & 0.25 & 0.30 & 0.35 & 0.40 & 0.45 & 0.50 \\
        \hline
20 & 0.0464 & 0.0519 & 0.0518 & 0.0478 & 0.0493 & 0.0582 & 0.0613 & 0.0698 & 0.0912 & 0.1071 & 0.1302 \\
50 & 0.0471 & 0.0476 & 0.0516 & 0.0512 & 0.0529 & 0.0666 & 0.0912 & 0.1239 & 0.1828 & 0.2613 & 0.3662 \\
100 & 0.0461 & 0.0512 & 0.0521 & 0.0581 & 0.0757 & 0.1078 & 0.1642 & 0.2552 & 0.3991 & 0.5643 & 0.7324 \\
    \end{tabular}
\end{table}

\begin{table}[H]
    \centering
        \caption{Empirical power of $T_2$ against the alternative  $\textrm{CN}_2((u,u),\rm{I}_2)$ ($n=20,50$ and $100$). }
    \label{tab:5.1.6}
    \begin{tabular}{r|rrrrrrrrrrr}
     $u$   &  0.00 & 0.05 & 0.10 & 0.15 & 0.20 & 0.25 & 0.30 & 0.35 & 0.40 & 0.45 & 0.50 \\
        \hline
20  & 0.0475 & 0.0519 & 0.0546 & 0.0482 & 0.0522 & 0.0567 & 0.0618 & 0.0724 & 0.0942 & 0.1078 & 0.1301 \\
50  & 0.0490 & 0.0495 & 0.0550 & 0.0523 & 0.0551 & 0.0693 & 0.0962 & 0.1277 & 0.1896 & 0.2664 & 0.3685 \\
100 & 0.0465 & 0.0519 & 0.0529 & 0.0587 & 0.0763 & 0.1076 & 0.1647 & 0.2567 & 0.3998 & 0.5634 & 0.7297\\
    \end{tabular}
\end{table}

\globalcolor{black}

\subsection{Discrete complex random variable}

Let
$$
Z = \left\{
\begin{array}{rl}
1+\ii, & \text{with probability } 1 /4, \\
1-\ii, & \text{with probability } 1 /4, \\
-1+\ii, & \text{with probability } 1 /4, \\
-1-\ii, & \text{with probability } 1 /4, \\
\end{array}
\right.
$$
be a complex random variable which is proper but not circularly symmetric (consider for instance $\e^{\ii\pi/4}Z$).

The empirical power for our test statistic $T_{n,\lambda}$ is given in Table~\ref{tab:table5.3.1}.

\begin{table}[H]
    \centering
        \caption{Empirical power of $T_{n,\lambda}$ for sample sizes $n=10,20$ and $50$, for $\lambda=1.0$, $0.1$ and $0.01$.}
    \label{tab:table5.3.1}
    \begin{tabular}{c|ccc}
    $n$ / $\lambda$ & 1.0 & 0.1 & 0.01 \\
    \hline
    10 & 0.233 & 0.057 & 0.053 \\
    20 & 0.526 & 0.057 & 0.055 \\
    50 & 1.000 & 0.056 & 0.053 \\
    \end{tabular}
\end{table}

The results of Table~\ref{tab:table5.3.1}  demonstrate the ability of our test to detect non-circularity for discrete complex random variables when the value of $\lambda$ is set to the default choice, $\lambda=1.0$
.

\textcolor{black}{We also consider the adjusted generalized likelihood
ratio test in \cite{OEK:2011} (see also \cite{OK:2009} and \cite{O:2008}), denoted GLRT, its robust version \cite{OKP:2011}, denoted RobGLRT,   and the Wald's type maximum likelihood test  of \cite{OEK:2011}, denoted WTML, for which empirical powers are presented in Table~\ref{tab:table5.3.2}.}

\globalcolor{black}

\begin{table}[H]
    \centering
        \caption{Empirical power of GLRT, RobGLRT and WTML for sample sizes $n=10,20$ and $50$.}
    \label{tab:table5.3.2}
    \begin{tabular}{c|ccc}
    $n$  & GLRT & RobGLRT & WTML \\
    \hline
    10 & 0.128 & 0.000 & 0.095 \\
    20 & 0.112 & 0.824 & 0.093 \\
    50 & 0.101 & 0.893 & 0.088 \\
    \end{tabular}
\end{table}

\globalcolor{black}

\subsection{A circularly symmetric random variable that does not have a density}

Consider the complex random variable $Z=\e^{\ii\Phi}$ with $\Phi\sim\rm{U}[-\pi,\pi)$. This non-Gaussian random variable {\it is} circularly symmetric but does not possess a density.

The empirical level  found for our test $T_{n,\lambda}$, with $\lambda=1.0, 0.1$ or $0.01$, is always between 0.05 and 0.057, for the sample sizes $n=10$, $20$ and $50$. These values are all very close to the nominal level $\alpha=0.05$. \textcolor{black}{For the other tests,  which do not rely on the bootstrap but on asymptotic distributions for the computation of their $p$-values, a larger sample size is necessary to attain their nominal level, as can be seen from Table~\ref{tab:table5.4}.}

\globalcolor{black}

\begin{table}[H]
    \centering
        \caption{Empirical levels for the tests  $T_{n,\lambda}$, GLRT, RobGLRT and WTML.}
    \label{tab:table5.4}
    \begin{tabular}{c|cccccc}
    $n$  & $T_{n,1.0}$ & $T_{n,0.1}$ & $T_{n,0.01}$ & GLRT & RobGLRT & WTML \\
    \hline
    10  & 0.0537 & 0.0501 & 0.0532 & 0.1179 & 0.0000  & 0.0853 \\
    20  & 0.0551 & 0.0561 & 0.0520 & 0.0753 & 0.0122 & 0.0626 \\
    50  & 0.0537 & 0.0550 & 0.0535 & 0.0601 & 0.0350 & 0.0544 \\
    100 & 0.0551 & 0.0528 & 0.0560 & 0.0543 & 0.0428 & 0.0526 \\
    \end{tabular}
\end{table}

\globalcolor{black}

\subsection{Contaminated distribution}

Consider the complex random variable
\begin{equation}\label{eq:sec4.4}
Z={\mathrm{P}}\e^{\ii\Theta}\text{ where }{\mathrm{P}}\sim\textrm{U}[0,1]\text{ and }\Theta\edist\left\{\begin{array}{rc}
0, & \text{with probability } 1/6, \\
2\pi/3, & \text{with probability } 1/6, \\
4\pi/3, & \text{with probability } 1/6, \\
2\pi U, &  \text{with probability } 1/2, \\
\end{array}\right.
\end{equation}
with $U\sim\textrm{U}[0,1]$ independently of $\textrm{P}$. This complex random variable is {\it not} circularly symmetric because it is contaminated, as clearly illustrated on Figure~\ref{fig:5.5.1}.

\begin{figure}[H]
    \centering
   \includegraphics[width=8cm]{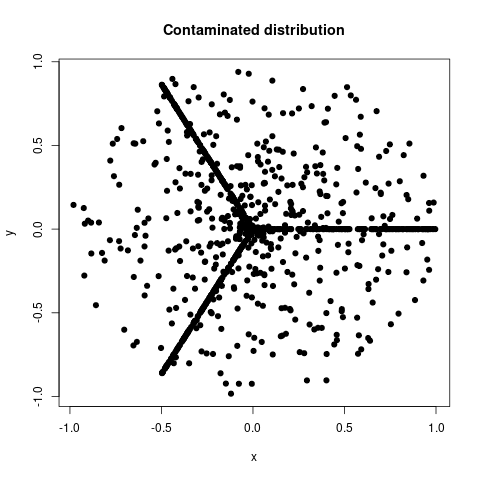}
    \caption{A sample of $n=1,000$ points generated according to the contaminated distribution~\eqref{eq:sec4.4}.}
    \label{fig:5.5.1}
\end{figure}

We randomly generate samples of sizes $n=10,20,50,100,200$ and $500$ and apply our test \textcolor{black}{as well as the GLRT, RobGLRT and WTML tests}. In Table~\ref{tab:5.5.1}, we see that our test (with $\lambda=1.0$) exhibits a power that increases with $n$, while the \textcolor{black}{other competing  tests have either no power or a power decreasing with $n$}. The other values $\lambda=0.1$ and $\lambda=0.01$ do not exhibit such a high power.

\globalcolor{black}

\begin{table}[H]
    \centering
        \caption{Power for the contaminated distribution~\eqref{eq:sec4.4}.}
    \label{tab:5.5.1}
    \begin{tabular}{c|cccccc}
         $n$   & 10    & 20    & 50    & 100   & 200 & 500 \\
         \hline
         $T_{n,1.0}$  & 0.0669 & 0.0657 & 0.0832 & 0.1078 & 0.2036 & 0.8318 \\
         $T_{n,0.1}$  & 0.0544 & 0.0541 & 0.0518 & 0.0563 & 0.0589 & 0.0620 \\
         $T_{n,0.01}$ & 0.0521 & 0.0561 & 0.0574 & 0.0536 & 0.0535 & 0.0580 \\
         GLRT         & 0.1889 & 0.1077 & 0.0698 & 0.0614 & 0.0504 & 0.0529 \\
         RobGLRT      & 0.0000 & 0.0623 & 0.0647 & 0.0578 & 0.0500 & 0.0518\\
         WTML         & 0.1217 & 0.0813 & 0.0600 & 0.0578 & 0.0484 & 0.0522
    \end{tabular}
\end{table}

\globalcolor{black}



\subsection{High dimensional complex random vector}

We consider $d$-dimensional complex normal vectors $Z\sim\textrm{CN}_d(0,\Gamma,P)$ as in~\cite{DLM:2016},  where we set $\Gamma$ as the $(d\times d)$-matrix that contains only ones, and where \textcolor{black}{$P=\ii A^\top A$, where the $d^2$ entries in the matrix $A$} are generated randomly (once for each value of $d$) from a $\mathcal{U}[0,1]$-distribution.

We then generated $M=1,000$ Monte-Carlo samples of observations from such random vectors, and considered the sample sizes $n=20,50,100$ and $200$  and the dimensions $d=2, 5, 10, 20, 50$ and $100$. We applied our test statistic $T_{n,\lambda}$ for the values $\lambda=1,0.1$ and $0.01$, which results in Tables~\ref{tab:highdim1}, \ref{tab:highdim01} and \ref{tab:highdim001}, respectively.

\globalcolor{black}

\begin{table}[H]
    \centering
        \caption{Empirical power based on $T_{n,\lambda}$, $\lambda=1.0$, for $d$-dimensional non-circular complex normal random vectors.}
    \label{tab:highdim1}
    \begin{tabular}{c|cccccc}
    $n$ / $d$ & 2 & 5 & 10 & 20 & 50 & 100 \\
    \hline
    20  & 0.268 & 0.187 & 0.074 & 0.065 & 0.039 & 0.049 \\
    50  & 0.752 & 0.516 & 0.081 & 0.063 & 0.046 & 0.053 \\
    100 & 0.993 & 0.923 & 0.078 & 0.063 & 0.049 & 0.054 \\
    200 & 1.000 & 1.000 & 0.129 & 0.066 & 0.056 & 0.063 \\
    \end{tabular}
\end{table}

\begin{table}[H]
    \centering
        \caption{Empirical power based on $T_{n,\lambda}$, $\lambda=0.1$, for $d$-dimensional non-circular complex normal random vectors.}
    \label{tab:highdim01}
    \begin{tabular}{c|cccccc}
    $n$ / $d$ & 2 & 5 & 10 & 20 & 50 & 100 \\
    \hline
    20  & 0.094 & 0.249 & 0.096 & 0.067 & 0.050 & 0.060 \\
    50  & 0.241 & 0.759 & 0.183 & 0.069 & 0.054 & 0.056 \\
    100 & 0.677 & 0.996 & 0.395 & 0.074 & 0.036 & 0.061 \\
    200 & 0.999 & 1.000 & 0.781 & 0.110 & 0.065 & 0.064 \\
    \end{tabular}
\end{table}

\begin{table}[H]
    \centering
        \caption{Empirical power based on $T_{n,\lambda}$, $\lambda=0.01$, for $d$-dimensional non-circular complex normal random vectors.}
    \label{tab:highdim001}
    \begin{tabular}{c|cccccc}
    $n$ / $d$ & 2 & 5 & 10 & 20 & 50 & 100 \\
    \hline
    20  & 0.067 & 0.08 & 0.075 & 0.083 & 0.066 & 0.049 \\
    50  & 0.079 & 0.117 & 0.091 & 0.08 & 0.053 & 0.051 \\
    100 & 0.095 & 0.229 & 0.183 & 0.102 & 0.056 & 0.065 \\
    200 & 0.113 & 0.769 & 0.529 & 0.209 & 0.055 & 0.055 \\
    \end{tabular}
\end{table}

\globalcolor{black}

As expected, the empirical power increases with $n$ and decreases with $d$.  \textcolor{black}{It also seems that smaller values of $\lambda$ might help the detection of non-circularity in higher dimensions to a certain degree.}


\medskip

\textcolor{black}{We close this section by revisiting the problem of choice of the weight parameter $\lambda$. While our results of Section ~\ref{sweight} certainly shed light on the intuition behind this parameter on the qualitative level, the practical problem of choosing $\lambda$ remains. To this end and in view of the disparity of results observed in our Monte Carlo study our suggestion is to try the test on a grid of values and choose a compromise value of $\lambda$ that renders the test powerful over a set of alternatives which are of potential interest. For instance, deviations from circularity within normality are better detected with a small value of $\lambda$, which might even be somewhat more robust to higher dimension, while for discrete alternative distributions or mixtures, a higher value of this parameter seems to be a better choice. Another idea is to implement a multiple test incorporating several values of $\lambda$ as in \cite{Tenreiro:11}. However, proper construction of such a multiple test requires a separate investigation.}

\section{Applications}\label{SEC6}

We consider raw METAR data consisting of wind direction (in degrees from North) and wind speed (in mph) for the first week of January 2020 in two Australian cities with different patterns of wind, namely the coastal city Sydney in New South Wales and the inland city Cloncurry in Queensland. The sample rate of Sydney records is roughly equal to two per hour (total sample size $n=360$), while for Cloncurry records it is around one per hour (total sample size $n=190$). No missing data were present.

We decided to store each observation as a complex number $z=x+\ii y=\rho \e^{\ii\vartheta}$, where $\rho$ is a measure of the wind speed (in mph) and $\vartheta$ is a measure of the wind direction (expressed in radians from East).

In order to apply our test, \textcolor{black}{as well as the  GLRT, WTML and RobGLRT tests}, to a reasonable sample size, and to remove potential outliers, we only kept low or moderate wind speeds. More specifically, we selected the subsamples of the above-mentioned two sets of data for which the Beaufort scale index is lower than or equal to 3 (i.e., up to a gentle breeze). For Sydney and Cloncurry, the \textcolor{black}{corresponding} cutoff value (wind speed $<$ 13) is roughly equal, respectively, to the empirical median and the third quartile of these two data sets.

\begin{figure}[H]
    \centering
    \includegraphics[width=8cm]{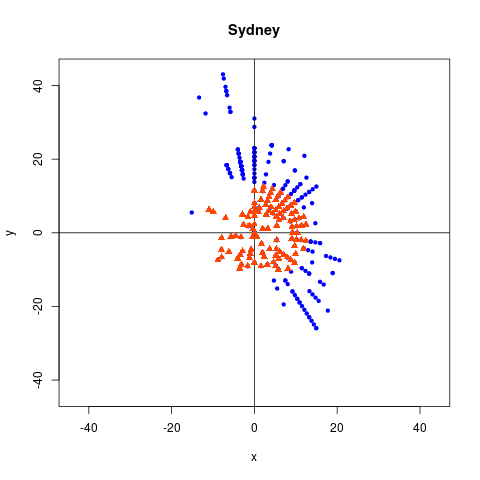}\includegraphics[width=8cm]{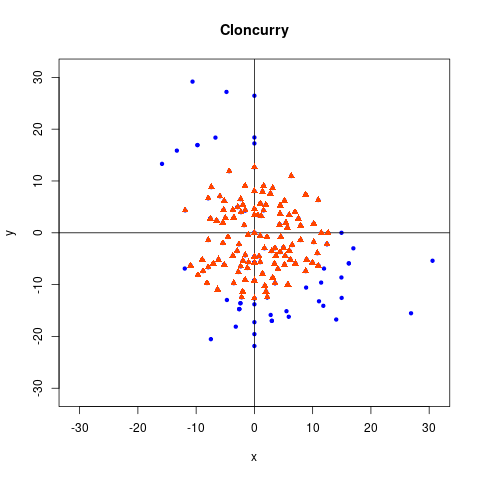}
    \caption{Wind data in Sydney and Cloncurry, first week of 2020 (source:  \url{https://mesonet.agron.iastate.edu/request/download.phtml?network=AU__ASOS}). Data points with a value of  $\sqrt{x^2+y^2}$ less (resp. greater) than 13 are displayed using orange triangles (resp. blue dots).}
    \label{fig:wind}
\end{figure}

We represent these data on the complex plane in Figure~\ref{fig:wind}, using orange triangles for the low-speed wind values as explained above. Both cities exhibit a clear noncircular pattern, more marked for Sydney than for Cloncurry.

\globalcolor{black}

We applied our test using our \texttt{R} package \texttt{CircSymTest}, as well as the tests GLRT, WTML and RobGLRT. All tests strongly reject circularity for Sydney, but only our test and the  RobGLRT test detect noncircularity for the (truncated) Cloncurry data set at the 5\% nominal level; see Table~\ref{tab:6.1}. Looking at the full data sets, noncircularity appears even more strikingly with very small $p$-values (not shown here) for all tests.

\begin{table}[H]
    \centering
        \caption{$p$-values when testing circularity for the truncated data sets.}
    \label{tab:6.1}
    \begin{tabular}{l|cccccc}
                           & $T_{n,0.01}$ & $T_{n,0.1}$ & $T_{n,1.0}$ & GLRT   & WTML   & RobGLRT \\
                           \hline
    Sydney ($n=178$)     &  0.000       &  0.000      &  0.000      &  0.006 &  0.007 &  0.011 \\
    Cloncurry ($n=148$)  &  0.001       &  0.003      &  0.019      &  0.258 &  0.260 &  0.014 \\
    \end{tabular}
\end{table}

\begin{figure}[H]
    \centering
    \includegraphics[height=8cm,height=14cm]{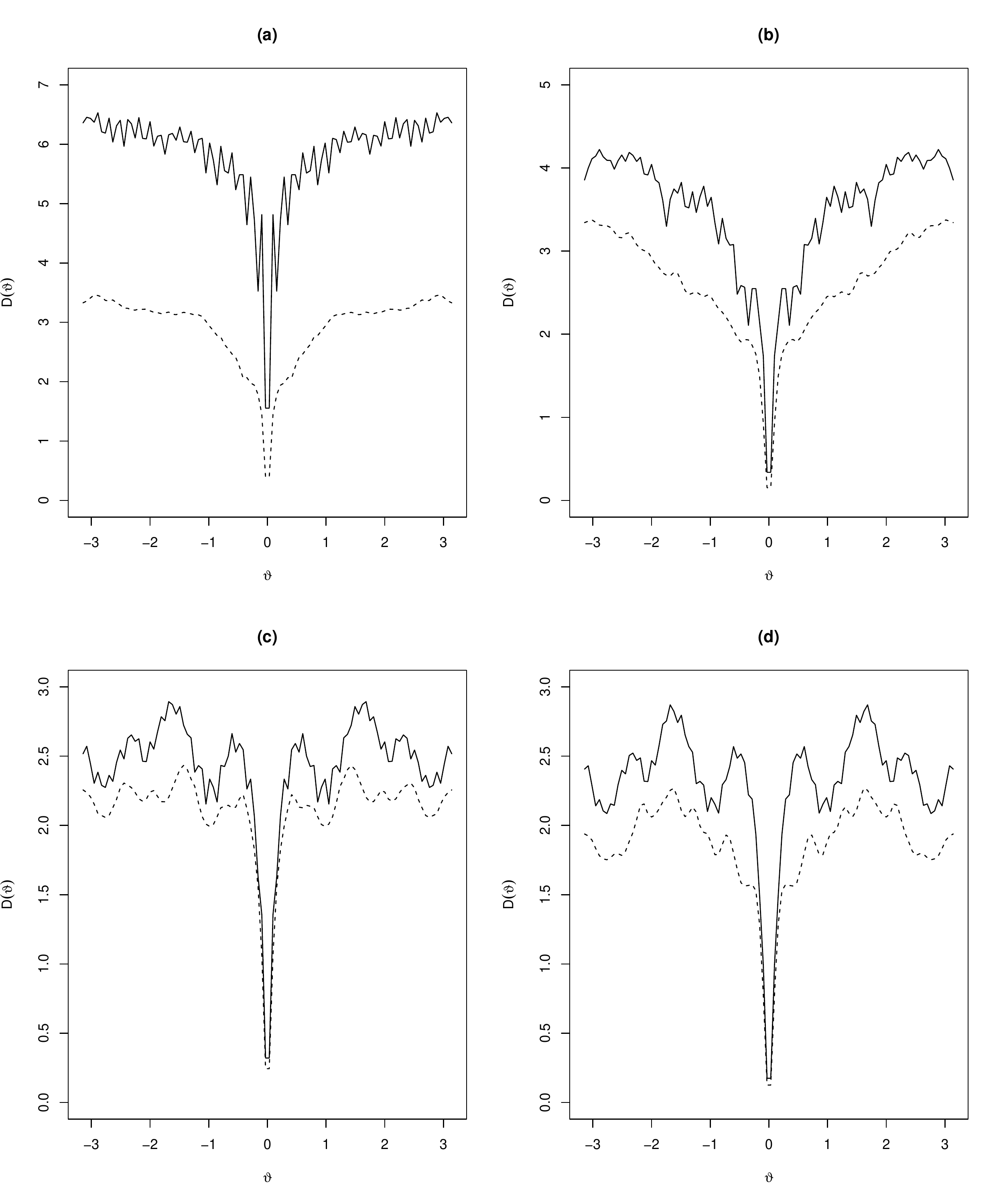}
    \caption{Value of $D(\vartheta)$ for $\vartheta\in[-\pi,\pi)$ at $\lambda=1.0$, both for the null (dashed-line) and non-null (solid line) hypotheses, for Sydney with $n=360$ (a) and $n=178$ (b), and for Cloncurry with $n=190$ (c) and $n=148$ (d).}
    \label{fig:6.2}
\end{figure}

It is even possible to investigate what values of $\vartheta$ in \eqref{NULL} are indicative of noncircularity by looking at the integral 
$D(\vartheta):=nT_{n,\lambda}(\vartheta)$, see  \eqref{ttheta}, for fixed value of $\lambda$. It is clear from  Figure~\ref{fig:6.2} that for $\lambda=1$, $D(\cdot)$ always takes largest values under the alternative hypothesis. For instance, for Cloncurry the largest discrepancies between the two curves are observed when $\vartheta$ takes values close to $\pi/6$ and $\pi/2$.

\globalcolor{black}

\textbf{Acknowledgements}: Research on this topic was initiated during the third  author's visit to the UNSW. Simos Meintanis would like to sincerely thank Pierre Lafaye de Micheaux and the School of Mathematics and Statistics of the UNSW for making this visit possible. \textcolor{black}{This research includes computations using the computational cluster Katana supported by Research Technology Services at UNSW Sydney.} \textcolor{black}{The authors would like to thank two anonymous referees for many valuable remarks.}



%
%
%

\begin{small}

\end{small}

\end{document}